\documentstyle[12pt]{article}
\newfont{\sm}{cmr10}
\def\a{\alpha}                        
\def\b{\beta}                         
\def\g{\gamma}
\def\l{\lambda}                       
                       \def\o{\omega}
                         
\def\O{\Omega}
\def\ho{\hat{\omega}}                 \def\ot{\otimes}
\def\hO{\hat{\Omega}}

\def\A{{\cal A}}              \def\MA{{\cal M}(\A)}

\newcommand{\bn}{\begin{equation}}     \newcommand{\ed}{\end{equation}}
\newcommand{\beq}{\begin{equation}}     \newcommand{\eeq}{\end{equation}}
\newtheorem{definition}{Definition}[section]
\newtheorem{proposition}{Proposition}[section]
\newtheorem{theorem}{Theorem}[section]

\newtheorem{corollary}{Corollary}[section]

\newcommand{\rf}[1]{(\ref{#1})}

\def\MlA{{\cal M}_\lambda({\cal A})}
\def\I{{ I}}
\def\C{{\bf C}}
\def\cXa{\overline{X}_\a}
\def\cXb{\overline{X}_\b}
\def\Xa{{\overline{X}_\a}}
\def\Xb{{\overline{X}_\b}}
\def\oX{\overline{X}}
\def\ver{{\rm ver}}
\def\codim{{\rm codim\ }}  
\def\tto{\Rightarrow}
\def\Ga{\Gamma_\A}

\def\ee{\emptyset}
\def\BA{{\cal B}_\A}
\def\VA{V_\A}

\def\Qui{Qui(\A)}  
\def\ba0{{\cal B}_\a^0} 
\def\qui{\Qui} 
\def\quila{ Qui_\l(\A)}
\def\la{\lambda}
\def\wE{\widetilde{E}}
\begin{document}
\begin{center}
{\Large\bf Non-resonance $D$-modules over arrangements of hyperplanes}
\end{center}
\bigskip
\begin{center}
Sergei Khoroshkin and Vadim Schechtman
\end{center}
\bigskip

\setcounter{footnote}{0}
\setcounter{equation}{0}

\section{Introduction}
\bigskip 

The aim of this note is a combinatorial description of a category 
of $D$-modules over an affine space, smooth along the stratification 
defined by an arrangement of hyperplanes. These $D$-modules are assumed 
to satisfy certain {\it non-resonance} condition. The main result, 
see Theorem 4.1, generalizes [3].

The work of S.K. was supported by the grants INTAS 93-10183 and 
RFBR 98-01-00303. 

\bigskip  

\section{Quivers of an arrangement}
\bigskip

{\large\bf 2.1.} Let $\A =\{H_i\}, i\in\I$  be a collection of affine  
hyperplanes in a complex
 affine space $X=\C^N$.  Hyperplanes $H_i$ define on $X$ the structure
 of a stratified space which we denote by the same letter $\A$ or by $X_\A$.
 The closure of a stratum $X_\a\subset \C^N$ is an intersection of some
 hyperplanes
 \bn
\cXa= \bigcap_{i\in J\subset I}H_i\ ,
\label{1.1}
\ed
and the stratum $X_\a$ is an interiour of $\cXa$:
$$X_\a=\cXa\setminus(\bigcup_{\cXb\subset \cXa, \cXb\not=\cXa}\cXb)\ .$$

The combinatorial structure of the adjacements of the strata can be described
 by a connected oriented graph $\Ga$. The vertices $\a$ of the graph correspond
 to the strata $X_\a$, the set of the vertices we denote by ${\rm ver}(\A)$;
 an arrow $\a\to\b \in {\rm arr}(\A)$ corresponds to the adjacement in
codimension $1$:
$$\a\to\b\qquad \Xa\supset\Xb\ , \quad \codim_\Xa\Xb=1\ .$$
Inclusions of closed strata define the partial order in the set
${\rm ver}(\A)$: $\a>\b \Leftrightarrow \Xa \supseteq \Xb
\Leftrightarrow$ there exists an oriented path from $\a$ to $\b$ in
 graph $\Ga$. This partial order is in agreement with integervalued
 function $\codim : \Ga \tto {\bf N}, \  \codim\a=\codim_{\C^N}\Xa$ which
 means that $\a >\b \tto \codim\a <\codim\b$. The partial order
 allows to attach to any vertice $\a$ of the graph $\Ga$ the full subgraph
$\Ga^\a\subset \Ga$ whose vertices are the vertices $\b\in \ver(\A)$
 which satisfy the condition $\b<\a$. One can easily see that the graph
 $\Ga^\a$ describes the combinatorial structure of affine subspace
 $\Xa$ with induced stratification.
 By $X_\ee$ we denote the the complement to all hyperplanes $H_i$ which is
 an open in $\C^N$ stratum.
\begin{definition}
A quiver $\VA=\{V_\a, A_{\a,\b}\}$ of the arrangement $\A$ is a collection of
 linear spaces $V_\a$ over $\C$, $\a\in\ver(\A)$ together with a collection of
 linear maps $A_{\a,\b}: V_\b\to V_\a,$ $\a,\b\in\ver(\A)$ such that
  $A_{\a,\b}=0$ unless $\a$ and $\b$ are  connected by an arrow
 ( $A_{\a,\b}\not = 0 \Rightarrow \a\to\b $ or $\b\to\a$) 
 , and which satisfy the following quadratic
 relations:
\bn
\sum_\b A_{\a,\b}A_{\b,\g}=0 \qquad {\mbox for\,  any}\; \a,\g\in \ver(\A),\;
 \a\not=\g\ .
\label{1.2}
\ed
\end{definition}
{\bf Remark. } One can see that the relation (\ref{1.2}) is nontrivial only in
two cases:

 {\it (i)} $\mid \codim\a-\codim\g\mid =2$ and $\a>\g$ or $\a<\g$;

 {\it (ii)} $\codim\a=\codim\g$, $\codim_{X_\a}X_\a\cap X_\g =1$.\\
In the last case there are one or two summands in the relation (\ref{1.2}).

The quivers $\VA$ compose the category $\Qui$ with natural morphisms. In other
 words, $\Qui$ is the category of finite-dimensional  modules over graded
 algebra $\BA$. The algebra $\BA$ is generated by orthogonal idempotents
$e_\a$, $\a \in \ver(\A)$ of degree 0, $\sum e_\a=1$ and by
degree one generators $a_{\a,\b}$, $\a,\b \in \ver(\A)$ such that
 $a_{\a,\b}=0$ if $\a$ and $\b$ are  connected by an arrow in $\Ga$ and,
 moreover,
$$ a_{\a,\b}e_\g= \delta_{\b,\g}a_{\a,\b},\qquad e_\g a_{\a,\b} =
 \delta_{\g,\a}a_{\a,\b}\ ,$$
$$\sum_\b a_{\a,\b}a_{\b,\g}=0,\qquad \a\not=\g.$$

Let $\VA=\{V_\a,A_{\a,\b}\}$ be a quiver of the arragement $\A$. We denote
 by $\VA^t=\{V_\a^t,A_{\a,\b}^t\}$ the following dual quiver:
\bn
V_\a^t=V_\a^*,\qquad
A_{\a,\b}^t=
(-1)^{(\codim\a+\codim\b-1)/2}A_{\b,\a}^*\ .
\label{1.3}
\ed
The correspondence $\VA\to \VA^t$ defines the contravariant functor of duality
$${}^t:\Qui\to \Qui\ .$$

%%%%%%%%%%%%%%%%%%%%%%%%%%%%%%%%
%%%%%%%%%%%%%%%%%%%%%%%%%%%%%%%

\bigskip

\setcounter{equation}{0}

\section{Quiver $D$-modules}
\bigskip

{\large\bf 3.1.\ }
For  any stratum $X_\a$ we denote by $\O_\a$ one-dimensional space
 of $\codim X_\a$-forms $\o_\a \in \det T^*_{X_\a}$, which are constant
 along $X_\a$, and by $\hat{\O}_\a$ denote one-dimensional space of top
 degree forms of the conormal bundle of $X_\a$, $\ho_\a \in 
T^*_{(X_\a\mid X)}$ constant along $X_\a$.

Let $\VA$ be a quiver of an arrangement $\A$. We relate to it a $D_X$-module
$E\VA$ by the following rules:

{\it (i)} $E\VA$ is generated by the linear space
 $\bigoplus_\a(V_\a\otimes\O_a)$ as $D_X$-module,

{\it (ii)} for any $\a$ and for any vectorfield $L_\a$ along $X_\a$ 
 with constant coefficients 
\bn
 L_\a(v_\a\ot \o_\a)=\sum_{\b:\a\to\b}
\frac{L_\a(\ho_\b)\wedge\o_\a}{\ho_\b\wedge\o_\b}
A_{\b, \a}v_\a\ot\o_\b ,
\qquad v_\a\in V_\a, \o_\g\in \O_\g, \ho_\g\in \hO_\g\ ,
\label{2.1}
\ed                                                    

{\it (iii)} for any $\a$ and for any affine function $f_\a$,
 $ f_\a\vert_{X_\a} =0$
\bn
f_\a(v_\a\ot w_\a)= \sum_{\b:\b\to\a} 
\frac{df_\a\wedge\ho_\b\wedge \o_\a}{\ho_\b\wedge \o_\b}
A_{\b, \a}v_\a\ot\o_\b,
\qquad v_\a\in V_\a, \o_\g\in \O_\g, \ho_\g\in \hO_\g\ ,
\label{2.2}
\ed                                                    
 In other words, $E\VA$ is a factor of free $D_X$-module 
$D_X\ot(\oplus_\a V_\a\ot\O_\a)$ modulo the relations \rf{2.1}--\rf{2.2}.

%%%%%%%%%%%%%       
%%%%%%%%%%
The closure $\Xa$ of any stratum is affine subspace of $X=\C^N$. We can
choose local affine coordinates $\{\overline{z}^\a,\overline{w}^\a\} =
\{ z_1^\a,\ldots , z^\a_{{\rm dim}X_\a},
w_1^\a,\ldots , w^\a_{{\rm codim}X_\a}\}$ such that $\Xa$ is defined by the
equations: 
\bn
w_1^\a=w_2^\a= \ldots = w^\a_{{\rm codim}X_\a}=0\ .
\label{2.3}
\ed
\begin{proposition}
{\em (i)} $E\VA$ is holonomic $D_X$-module with regular singularities
 flat along the stratification $\A$ (that is,
$ ss.E\VA \subset \bigcup_\a T^*_{(X_a\mid X)}$);

{\em (ii)} The global sections of $E\VA$ are:
\bn
\Gamma (E\VA)\equiv
\oplus_\a\C[\{\overline{z}^\a\},
\{\frac{\partial}{\overline{w}^\a}\}] V_\a\ot\O_\a\ .
\label{2.4}
\ed
\end{proposition}
{\bf Proof .} {\it (ii)} Due to the relaions \rf{2.1} and \rf{2.2}
the space of global sections is isomorphic to a factorspae of linear
 space
$$U=\oplus_\a U^\a=\oplus_\a\C[\{\overline{z}^\a\},
\{\frac{\partial}{\overline{w}^\a}\}] V_\a\ot\O_\a \ .
$$                                           
So it is sufficient to define the structure of $D_X$ module in
 linear space $U$. We can equip $U$ with a structure of graded space
by the relations ${\rm deg\ }(v_\a\ot w_a) =0, {\rm deg\ }(z^\a_i) =
{\rm deg\ }(\partial_{w_j^\a}) =1$ $(\partial_{w_j^\a}=
\frac{\partial}{\overline{w}^\a})$ and define an action of the
 generators $x_i, \partial_{x_i}$ of the ring $D_X$ by induction
over the grading. The generators $\partial_{z_i^\a}$ and $w_j^\a$
 act on the elements of degree 0 $v_\a\ot \o_\a$ according to the relations
\rf{2.1} and \rf{2.2}, whereas the generators
 $\partial_{w_j^\a}$ and $z_i^\a$ act on them by free multiplication.
 If we already define the action
 of the generators $x_i, \partial_{x_i}$ on the elements of degree
 $\leq (n-1)$, the action of $x_i, \partial_{x_i}$ on the elements
 of degree $n$ from $U_\a$ is uniquely defined by the following conditions:\\
 {\it (a)} the elements
 $\partial_{w_j^\a}$ and $z_i^\a$ act on subspace $U^\a$ by free
 multiplication, in particular,
 \bn
 [z_i^\a,z_k^\a]=[z_i^\a,\partial_{w_j^\a}]=
 [\partial_{w_j^\a},\partial_{w_l^\a}] =0\ ,
 \label{2.5}
 \ed
 {\it (b)} the operators $\partial_{z_i^\a}$ and $w_j^\a$   commute with
 $\partial_{w_j^\a}$ and $z_i^\a$   according to the relations in $D_X$,
 that is,
 \bn
 [\partial_{z_i^\a},\partial_{w_j^\a}] =[w_j^\a,z_i^\a]=0\ ,
 [\partial_{z_i^\a},z_k^\a]=\delta_{i,k}, \
 [w_j^\a, \partial_{w_l^\a}]=-\delta_{j,l}\ .
 \label{2.6}
 \ed
 We have to prove that this action satisfied all the relation in $D_X$,
 so we are to check the relations
 \bn
 [\partial_{z_i^\a},\partial_{z_k^\a}] =[w_j^\a,w_l^\a]=
 [\partial_{z_i^\a},w_j^\a]=0\ .
 \label{2.7}
 \ed
 The relations \rf{2.5} and \rf{2.6} allow to reduce the check of
 \rf{2.7} to the space of degree zero, where it is valid due to
 \rf{2.1}, \rf{2.2} and \rf{1.1} (this is a simple but crucial calculation).

{\bf (i)} The description of the global sections of the module $E\VA$ enables
 one to describe the concrete good filtration of $E\VA$. Remind that a
filtration $M_k$ of $D_X$-module $X$ is good, if it is an argeement with
 a filtration  $D_k$ of $D_X$ over the degree of differential operator:
 $D_kM_l\subset M_{k+l}$; $D_kM_l=M_{k+l}$ for all $l\geq l_0, k\geq 0$;
 and all $M_k$ are finite $O_X$-modules.

Let, as before$ U^\a =\oplus_\a\C[\{\overline{z}^\a\},
\{\frac{\partial}{\overline{w}^\a}\}] V_\a\ot\O_\a \ .
$. We showed that, globally, $\Gamma(E\VA)=U=\oplus_\a U_\a$. We 
 equip $U_\a$ with a filtration putting ${\rm deg}(V_\a\ot\O_\a)=
\codim\a\ ,{\rm deg\ }z_i^\a=0\ , {\rm deg\ }\partial_{\o_i^\a} = 1 ,
 U^\a_k = \{ m\in U_\a , {\rm deg\ }m \leq k\} $ and let $U_k = 
 \oplus_\a U^\a_k$. Then, due to the relations \rf{2.1} and \rf{2.2},
 the filtration $0\subset U_0 \subset \ldots \subset U$ is a good filtration.
 Moreover, the same relations show that the product
$$\prod_\a\left(\prod_i p_i^\a\cdot\prod_j w_j^\a\right)$$
annihilate ${\rm gr\ }E\VA$ (here $p^\a_i ={\rm gr\ }\frac{\partial}
{\partial_{z_i^\a}}$). It means that $E\VA$ is holonomic $D_X$ module
 with singular support contained in $\cup_\a T^*_{(X_\a\vert\C^N)}$.
 It has regular singularities due to \rf{2.1}, \rf{2.2}                    
                     
Thus we have a functor
\bn E:\Qui\longrightarrow \MA\ ,
\label{E}
\ed
where $\MA$ is a category of holonomic RS $D_X$-modules flat anong
 the stratification $\A$.  Note that the presentation of the module $E\VA$
 allows one in addition easily describe the sections over the open sets
 $$X_H=X\setminus H$$
 where $H$ is a hyperplane.
Let, for instance, $H_\b$ be a generic hyperplane which contains the stratum
 $X_\b$. Then the global sections over $X_{H_\b}$ constitue the linear space 
 $$U_\b=\oplus_{\a: \b\not\to\a}U^\a=\oplus_\a\C[\{\overline{z}^\a\},
\{\frac{\partial}{\overline{w}^\a}\}] V_\a\ot\O_\a \ ,$$
 where the polynomial over $z_i^\a$ should be replaced by the regular
functions of $\Xb\cap U$.
%%%%%%%%%%%%%%%%%%%%%%%%%%%%%%%%%%%%%%%%%%%%%
%%%%%%%%%%%%%%%%%%%%%%%%%%%%%%%%%%%%%%%%%%%
\smallskip

{\large\bf 3.2.} One can also localize the whole construction to the
 complement of certain hyperlanes, that is, take as an initial space $X$
 the complement to some hyperplanes instead of affine space.
For instance, the configuration $\{\A, H\}$ 
of hyperplanes $\A=\{ H_i\}, i\in  I$   in $\C^N\setminus H$ defines the 
 quivers which consist of linear spaces $V_\a, \a\in \Ga$,  
 together with
 linear maps $A_{\a,\b}:V_\b \to V_\a$, which are zero unless
 $\a$ and $\b$ are  
 connected by an arrow, and linear maps $N_\a:V_\a\to V_\a$ which are not zero
 only for the strata $X_\a$ which intersect with $H$ over codimension one.
  These operators are subjected to the relations:
\bn
\sum_\b A_{\a,\b}A_{\b,\g}=0 \qquad {\mbox for\,  any}\; \a,\g\in \ver(\A),\;
 \a\not=\g\ .
\label{3.1}
\ed
\bn
N_\b A_{\b,\a}=A_{\b,\a}\left(N_\a+\sum_{{\g\not =\b:}\atop
{\overline{X}_\g\cap H=\Xb\cap H}}A_{\a,\g}A_{\g,\a}\right)\ ,\qquad 
\a\to\b\ ,
\label{3.2}
\ed 
\bn
 A_{\a,\b}N_\b=\left(N_\a+\sum_{{\g\not =\b:}\atop
{\overline{X}_\g\cap H=\Xb\cap H}}A_{\a,\g}A_{\g,\a}\right)
A_{\a,\b}\ ,\qquad 
\a\to\b\ ,
\label{3.3}
\ed
The corresponding $D_X$-module is described by analogous relations:
$$
 L_\a(v_\a\ot \o_\a)=\sum_{\b:\a\to\b}
\frac{L_\a(\ho_\b)\wedge\o_\a}{\ho_\b\wedge\o_\b}
A_{\b, \a}v_\a\ot\o_\b +$$
\bn +\phi_{L_\a}(v_\a\ot\o_\a),
\qquad v_\a\in V_\a, \o_\g\in \O_\g, \ho_\g\in \hO_\g\ ,
\label{3.4}
\ed                                                    
for any $\a$ and for any vectorfield $L_\a$ along $X_\a$ 
 with constant coefficients
\bn
f_\a(v_\a\ot w_\a)= \sum_{\b:\b\to\a} 
\frac{df_\a\wedge\ho_\b\wedge \o_\a}{\ho_\b\wedge \o_\b}
A_{\b, \a}v_\a\ot\o_\b,
\qquad v_\a\in V_\a, \o_\g\in \O_\g, \ho_\g\in \hO_\g\ ,
\label{3.5}
\ed                                                    
  for any $\a$ and for any affine function $f_\a$,
 $ f_\a\vert_{X_\a} =0$.

The summand $\phi_{L_\a}(v_\a\ot\o_\a)$ is nonzero only for those $\Xa$
 which intersect with $H$ over in the codimension $1$ and looks like
$$
\phi_{L_\a}(v_\a\ot\o_\a)=
\frac{L_\a(f)}{f}\left(N_\a v_\a\ot \o_\a-
\right.
$$
 $$ \left.\sum_{{\g:\overline{X}_\g\cap H=\Xa\cap H}}
\frac{df\wedge \ho_\g\wedge\o_{H\a}}{\ho_\gamma\wedge\o_\gamma}\cdot
\frac{\ho_\a\wedge\o_\a}{df\wedge\ho_\a\wedge\o_{H\a}}A_{\gamma,\delta}
A_{\delta,\a}v_\a\ot\o_\gamma\right)
$$
Here $\o_{H\a}$ is an arbitrary top form on $H\cap X_\a$ with constant 
coefficients, the symbol $\delta$ stands for a stratum which coincides
 with $X_\a+X_\g$; if $X_\a+X_\g$ is not a stratum then the corresponding 
 summand in r.h.s. of the above formula is zero ($X_\a+X_\g$
means a sum of linear subspaces in linear space with an origin in 
$X_\a\cap X_\g$). The statement of Proposition remains true for such 
 $D$-modules.
%%%%%%%%%%%%%%%%%%%%%%%%%%%%%%%%%%%%%%%%%%%%%%%%%%
%%%%%%%%%%%%%%%%%%%%%%%%%%%%%%%%%%%%%%%%%%%%%%%%%
\bigskip

{\large\bf 3.3. A free resolution of quiver $D$-modules.}

\bigskip

Let $T_\a=T'_\a\oplus T''_\a$, $dim\ T_\a=N$ be a direct sum of the space of 
 the sections of tangent bundle $T_{X_\a}$ to a stratum $X_\a$ 
constant along $X_\a$
 and of the space of sections of conormal bundle $T^*_{(X_\a|{\bf C}^N)}$
 constant along $X_\a$.
In local coordinates $\overline{z}^\a,\overline{w}^\a$ \rf{2.3} the spaces
 $T'_\a$ and $T''_\a$ are: $T'_\a=\C<\partial_{z_i^\a}>$, 
 $T''_\a=\C<d{w_j^\a}>$. In other way the space $T_\a$ can be identify
 with the space of constant sections of tangent bundle $T_{T^*_{(X_\a|\C^N)}}$
to the space of conormal bundle $T^*_{(X_\a|\C^N)} \subset T^*_{\C^N}$.

There is a natural embedding $I_\a$ of the space $T_\a$ into the space 
 of  first order differential operators:
$$I_\a:T_\a\hookrightarrow D_{\leq 1}\ .$$
In local coordinates \rf{2.3} $I_\a(\partial_{z_i^\a} = \frac{\partial}
{\partial z_i^\a}$, $I_\a(dw_j^\a = w_j^\a$. The map $I_\a$ does not depend 
 on ths choise of affine coordinates $\overline{z}^\a,\overline{w}^\a$
 satisfying \rf{2.3}.

{} For any quiver $\VA$ the relations \rf{2.1} and \rf{2.2} define 
also a natural map 
$$Q_\a:\ T_\a\ot V_\a\ot\Omega_\a\to\oplus_\b V_\b\ot\Omega_\b$$
which looks as follows:
\beq
Q_\a(\xi\ot v_\a\ot\o_\a)=\sum_{\b:\a\to\b}\frac{\xi(\ho_\b)\wedge\o_\a}
{\ho_\b\wedge\o_\b}A_{\b,\a}v_\a\ot\o_\b\ ,\qquad
\xi\in T'_\a, \ho_\b\in\hO_\b\ ,
\label{2.8}
\eeq
\beq
Q_\a(df_\a\ot v_\a\ot\o_\a)=\sum_{\b:\b\to\a}
\frac{df_\a\wedge\ho_\b\wedge\o_\a}
{\ho_\b\wedge\o_\b}A_{\b,\a}v_\a\ot\o_\b\ ,\qquad
df_\a\in T''_\a, \ho_\b\in\hO_\b\ .
\label{2.9}
\eeq
Let $E^n\VA, n=1,2,\ldots , N$ be a free left $D$-module
$$E^n\VA=\oplus_\a D\ot\left(\bigwedge^n(T_\a)\ot V_\a\ot\O_\a\right)$$
and $d^{(n)}:\ E^n\VA\to E^{n-1}\VA$ be the following 
$D$-linear map of free $D$-modules:
$$
d^{(n)}(1\ot t_1\wedge\ldots \wedge t_n\ot v_\a\ot \o_\a)=
\sum_{k=1}^n (-1)^{k+1}I_\a(t_k)\ot t_1\wedge\ldots\wedge\widehat{t_k}
\wedge\ldots\wedge t_n\ot v_\a\ot\o_\a-$$
\beq
-\sum_{k=1}^n (-1)^{k+1}1\ot t_1\wedge\ldots\wedge\widehat{t_k}
\wedge\ldots\wedge t_n\ot Q_\a(t_k\ot v_\a\ot\o_\a) \ .
\label{2.10}
\eeq
We state that $d^{(n-1)}d^{(n)} =0$ and
 \begin{proposition}
$\left(E^{\bullet}\VA,d^{({\bullet})}\right)$ is a free resolution of 
$D$-module $E\VA$.
\end{proposition}
{\bf Proof}. The relation $d^{(n-1)}d^{(n)} =0$ is equivalent to basic
 relations \rf{1.2} on the operators $A_{\a,\b}$. The proof of the exactness
 of the complex is the same as for standard Koszul resolution of the ring
 ${\cal O}_X$.
\begin{proposition}
The dual complex 
$$\bigwedge^N(T^*_{\C^N})[-N]\ot_{{\cal O}}{\rm Hom}_D
\left(E^{\bullet}\VA,d^{({\bullet})}\right)$$
is a free resolution of the module $E\VA^t$.
\end{proposition}

{\bf Proof}.
One can see that 
\bn
{\rm Hom}_D(D\ot\bigwedge^n(T_\a)\ot V_\a\ot\O_\a,D)\approx
D\ot\bigwedge^{N-n}(T_\a)\ot V_\a^*\ot\O_\a\ot\bigwedge^N(T_{\C^N})
\label{2.11}
\ed
and the dual maps $d^{(n)*}$ of right $D$-modules 
$D\ot\bigwedge^{N-n}(T_\a)\ot V_\a^*\ot\O_\a\ot\bigwedge^N(T_{\C^N})$
are of the form analogous to \rf{2.10}:
$$
d^{(n)*}(1\ot t_1\wedge\ldots \wedge t_n\ot v_\a^*\ot \o_\a\ot \pi_\ee)=
\sum_{k=1}^n (-1)^{k+1}I_\a(t_k)\ot t_1\wedge\ldots\wedge\widehat{t_k}
\wedge\ldots\wedge t_n\ot v_\a^*\ot\o_\a\ot\pi_\ee-$$
 $$
-\sum_{k=1}^n (-1)^{k+1}1\ot t_1\wedge\ldots\wedge\widehat{t_k}
\wedge\ldots\wedge t_n\ot Q_\a^*(t_k\ot v_\a^*\ot\o_\a)\ot\pi_\ee \ ,
 $$
where $Q_\a^*$ is given by the relations \rf{2.8}, \rf{2.9} with
 $v_\a^*, v_\b^*$ and $A^*_{\b,\a}$ substituted instead of
$v_\a, v_\b$ and $A_{\a,\b}$, $\pi_\ee$ is an $N$-vector on $\C^N$.
 The passage from the right $D$-modules to the left $D$-modules
 gives the sign $(-1)^{(\codim \a-\codim \b -1)/2}$ before
 $A^*_{\b,\a}$.
\begin{corollary}
The $D$-module $E\VA^t$ is the dual to $E\VA$. 
\end{corollary}
  Note that the appearence of $N$-vector in \rf{2.11} instead of the 
 top form shows that it is more natural to treate the quiver's 
 $D$-modules as right $D$-modules but the presense of affine structure
 in the stratified space allows one to pass from the left to right
 $D$-modules just by simple change of signes in their defining relations.

%%%%%%%%%%%%%%%%%%%%%%%%%%%%%%%%%%%%%%%%%%%%%%%%%%%%%%%%%
%%%%%%%%%%%%%%%%%%%%%%%%%%%%%%%%%%%%%%%%%%%%%%%%%%%%%%%

\bigskip

\setcounter{equation}{0}

\section{Weighted quivers and weighted $D$-modules}
\bigskip

{\large\bf 4.1. The Category $\quila$ of weighted quivers.}
\bigskip 

 Let us assign to any hyperplane $H_i$ of the arrangement $\A$
 a complex number (weight) $\la_i$, and set
\beq
\la_\a=\sum_{i:H_i\supset \Xa}\la_i
\label{2.2.1}
\eeq
 for any stratum $X_\a$ of the arrangement, and
\beq
\la_{\a,\b}=\sum_{i:H_i\supset \Xb, H_i\not\supset X_\a}\la_i
\label{2.2.2}
\eeq
for any arrow $\a\to\b$ of the graph $\Gamma_\A$.
Let us fix a stratum $\Xa$. Then we can attach to this stratum a
 subalgebra  $\ba0\subset {\cal B}_\A$  of the path algebra
 ${\cal B}_\A$ (see section 1) generated by the elements $A_\a^\b =
A_{\a,\b}A_{\b,\a}$, where $\b$ runs  all the vertices of $\Ga$ 
satisfying the condition $\a\to\b$. The algebra $\ba0$ can be described
 as an algebra with generators $A_\a^\b, \a\to\b$ subject to 
 quadratic relations 
\beq
[\sum_{\b:\a\to\b\to\g}A_\a^\b,A_\a^\delta]=0
\label{2.2.3}
\eeq
which are enumerated by the flags $\a\to\delta\to\gamma$.
One can easily see that for any
quiver $V_\A$ the linear space $V_\a$ is equiped with a structure of 
 $\ba0$ module. We define a subcategory $\quila\subset \qui$ as follows:

\begin{definition}
A quiver $\VA \in \Qui$ belongs to the subcategory $\quila$ if for any stratum
$X_\a$, the linear space $V_\a$ considered as a $\ba0$-module admits a 
 finite filtration
 by $\ba0$-submodules with the factors isomorphic to 
${\bf C}_{\l,\a}$, where ${\bf C}_{\l,\a}$ is one-dimensional $\ba0$-module
 in which the generators $A_\a^\b$ act by multiplication over $\l_{\a,\b}$:
\beq
\rho_{\l,\a}(A_\a^\b)=\l_{\a,\b}\cdot{\rm Id} .
\label{2.2.4}
\eeq
\end{definition}

In other words, the action of the operators $A_{\a,\b}$ in the space $V_\a$
 can be given in some basis by triangular matrices with diagonal
 $\l_{\a,\b}\cdot{\rm Id}$. The proposition below describes the properties
 of the category $\quila$.

\begin{proposition}
{\em (i)} The category $\quila$ is a full abelian subcategory of $\Qui$
 closed under extensions;\\
{\em (ii)} The category $\quila$ is stable with respect to the  
involution ${}^t$;\\
{\em (iii)} The category $\quila$ has finite number of irreducibles; the 
 irreducibles $L_{\l,\a}$ are labeled by the vertices $\a$ of the graph $\Ga$
 and are supported on the subgraph  $\Ga^\a$ of the closure of the stratum
 $X_\a$, that is $(L_{\l,\a})_\b\not= 0\Rightarrow \a>\b$ and 
$(L_{\l,\a})_\a \not= 0$.
\end{proposition}

 The conjecture that the category $\quila$ containes the projective covers
 and has homological dimension at most $N$ also makes sense but we did not
 try to prove it.
\bigskip

{\bf Proof}. The statements {\it (i)} and {\it (ii)} follow directly 
from the definitions and from the description of the involution. The most 
 nontrivial statement {\it (iii)} is based on the explicit construction of
 an  analog of Verma modules $M_{\l,\a}\subset \quila$
 which we give  below.
The quiver $M_\l=M_{\l,\ee}$ has the unique irreducible 
factor $L_\l = L_{\l,\emptyset}$
 and is characterized by the property
\beq
{\rm Hom}_{\Qui}\left(M_\l,\VA\right)\approx{\rm Hom}_{{\rm B}^0_\ee}
\left({\bf C}_{\l,\ee},
(\VA)_\ee\right)
\label{2.2.5}
\eeq
for any quiver $\VA\in\Qui$. The ${\rm B}^0_\ee$ 
 character of ${\bf C}_{\l,\ee}$ is defined by
the relation \rf{2.2.4}:
\beq
\rho_{\l,\ee}(A_{\ee,i}A_{i,\ee})=\l_i\cdot{\rm Id}\ .
\label{2.2.6}
\eeq
The construction of the quivers $M_{\l,\a}$ and of their irreducible factors
 uses the inclusions $j_\a:\Xa\hookrightarrow{\bf C}^N$. The inclusion $j_\a$
 defines the induced stratification $\A^\a$ in affine space $\Xa$ which
 is described by the graph $\Ga^\a$. It also defines the full exact functor
 (which we denote by the same letter)
$$j_\a:\ Qui(\A^\a)\to\Qui\ .$$
A collection of weights $\l=\{\l_i\}$ for the stratification $\A$ naturally
 defines a collection of weights $\l^\a=\{\l_{\a,\b}\}$ where $\b$ runs
 all the vertices of $\Ga^\a$ satisfying the condition $\a\to\b$. Simple
 linear algebra calculations show that the functor $j_\a$ maps the category
 $Qui_{\l^\a}(\A^\a)$ into $\quila$:
\beq
j_\a\left(Qui_{\l^\a}(\A^\a)\right)\subset \quila\ .
\label{2.2.7}
\eeq
Thus the irreducible quiver $L_{\l,\a}$ with a support $\Ga^\a$ is
$$L_{\l,\a}=j_\a(L_{\l^\a})$$
and the proof of the statement {\it (iii)} of the proposition follows 
by induction over the codimension of the support of the irreducibles
 using \rf{2.2.5}. Let us describe now the quiver $M_\l$.

\bigskip

{\bf Construction of the quiver $M_\l$}.
\bigskip

 The quiver $M_\l$ which represents the functor 
$${\cal F}(\VA)={\rm Hom}_{{\rm B}^0_\ee}
\left({\bf C}_{\l,\ee},(\VA)_\ee\right)$$
can be defined in terms of the flag complex ~\cite{SV}.
\begin{definition}
A (complete) flag $L_{\a_0,\ldots ,\a_n}$ in the stratified space 
$X_\A$ is a chain
 of the closures of strata
\bn
\oX_{\a_0}\supset\oX_{\a_1}\supset\ldots\supset\oX_{\a_n}\ ,\qquad
{\rm codim}\oX_{\a_j}=j\ .
\label{2.2.8}
\ed
\end{definition}

In terms of the graph $\Ga$, the complete flag $L_{\a_0,\ldots ,\a_n}$
 is a chain $\a_0\to \ldots\to \a_n$, $\a_0=\ee$.  
The linear space $(M_\l)_\g$ of the quiver $M_\l$ we identify with the 
 factorspace of  linear space with a basis given by all flags
$L_{\a_0,\ldots ,\a_n}$, $\a_n=\g$ modulo the relations
\bn
\sum_{\b:\a_{k-1}\to\b\to\a_{k+1}}
L_{\a_0,\ldots ,\a_{k-1},\b,\a_{k+1},\ldots\a_n}=0
\label{2.2.9}
\ed
valid for all partial flags 
$\a_0\to\ldots \to\a_{k-1}>
\a_{k+1}\to\ldots\to\a_n$, $\a_0=\ee, \a_n=\g, 
\codim \a_{k+1}= \codim \a_{k-1}+2$.
 For any arrow $\a\to\b$ the linear operators 
$$A_{\b,\a}:(M_\l)_\a\to(M_\l)_\b,\qquad
A_{\a,\b}:(M_\l)_\b\to(M_\l)_\a$$
 have the following form:
\bn
A_{\b,\a}(L_{\a_0,\ldots,\a_{n-1},\a})=
L_{\a_0,\ldots,\a_{n-1},\a,\b}\ ,
\label{2.2.10}
\ed
$$
A_{\a,\b}L_{\a_0,\ldots,\a_{n},\b}=
$$
\bn
\sum_{k=0}^n(-1)^k\l_{\a_{n-k},\a_{n-k+1}}
\sum_{\a'_{n-k+1},\ldots,\a'_n:\a'_j\to\a_{j+1},\a'_j\not=\a_j,\a'_n=\b}
L_{\a_0,\ldots,\a_{n-k},\a'_{n-k+1},\ldots,\a'_n}\ .
\label{2.2.11}
\ed
One can verify by linear algebra calculations that the operators \rf{2.2.10}
 and \rf{2.2.11} are correctly determined maps of corresponding factorspaces
 and define a quiver from $\quila$ satisfying the condition \rf{2.2.5}.
%%%%%%%%%%%%%%%%%%%%%%%%%%%%%%%%%%%%%%%%%%%%%%%%%
%%%%%%%%%%%%%%%%%%%%%%%%%%%%%%%%%%%%%%%%%%%%%%%%%

\bigskip

{\large\bf 4.2. The category $\MlA$ of weighted $D$-modules.}

\bigskip

{} For any collection $\l=\{ \l_i, i\in I\}$ of weights we define also
 a category $\MlA\subset {\cal M}(\A)$ of weighted $D$-modules.

Let $f_\g$, $\g\in\Ga$ stands for a generic linear (affine) function over
 $\C^N$, such that $f_\g\vert_{\oX_\g}= 0$. Let $\Phi_{f_\a}(M)\vert_{X_\a}$ be
 the restriction to the open part $X_\a$ of affine subspace $\oX_\a$ of
 $D$-module of vanishing cycles $\Phi_{f_\a}(M)$. If $M\in {\cal M}(\A)$ then 
 $\Phi_{f_\a}(M)\vert_{X_\a}$ is a flat $D_{X_\a}$-module which one can
 describe by a flat connection in trivial bundle over $X_\a$ with 
regular singularities on $\oX_\a$. Let us fix for a moment a collection of
 generic functions $f_\a, \a\in \Ga$.

\begin{definition}
A $D$-module $M$ from ${\cal M}(\A)$ belongs to the category $\MlA$, if for any
 stratum $X_\a$ $D_{X_\a}$-module $\Phi_{f_\a}(M)\vert_{X_\a}$ admits  
 $D_{X_\a}$-filtration with factors isomorphic to the modules given by a
 flat connection in one-dimensional bundle with a form
$$\nu_\a=\sum_{\b:\a\to\b}\l_{\a,\b} d\log f_\b$$
\end{definition}

In other words, for some affine coordinates $D$-module 
$\Phi_{f_\a}(M)\vert_{X_\a}$ is given by a connection in trivial bundle with
 a form 
$$\nu_\a=\sum_{\b:\a\to\b}L_{\a,\b} d\log f_\b$$
where $L_{\a,\b}$ is a triangular matrix with constant coefficients 
 and diagonal entries
 being equal to $\l_{\a,\b}$.
 One can see that the definition of the category $\MlA$ does not depend 
 on the choise of generic linear functions $f_\a$, $\a\in\Ga$ and depends
 only on the projection $\tilde{\pi}(\l)$ of the collections of weights to
 $\C/{\bf Z}$, $\tilde{\pi}(\l)=\exp(2\pi i\l)$.

\begin{definition}
A collection $\{\l\}$ of weights is called {\bf  non-resonance} if for
 any $\a, \b, \g, \delta \in \Ga$, the following conditions hold:
\bn
\l_{\a,\b}\not\in {\bf Z}\setminus \{ 0\}\qquad {\rm and} \qquad
\l_{\a,\b}-\l_{\g,\delta}\not\in{\bf Z}\ .
\label{nonresonanse}
\ed
\end{definition}

The non-resonance property also essentially
 depends only on the exponents of weights in a sense that one can
 or cannot get the condition \rf{nonresonanse} by shifting some of $\l_i$
 by an integer. 
One of the main properties of weighted $D$-modules is given in 
 the following theorem.
  
\begin{theorem}
{} For any  non-resonance collection $\{\l\}$ of weights, the functor
 $E$ establishes an equivalence of categories 
  $\quila$ and $\MlA$.
\end{theorem}

The statement of the theorem is a direct generalization of the results
 of \cite{Kh}. Let us describe the scheme of the proof. It consists of the
inductive (over the codimension of strata) construction of an inverse
 functor $\wE : \MlA\to\quila$ with simultenious proof of basic statement
 that the functor $E$ coincides with Beilinson's functor \cite{Be}
 of glueing
the  $D$-modules. The functor $\wE$ is as follows. The linear space 
 $(\wE M)_\a$ for $M\in \MlA$ is isomorphic to the space of flat sections of  
 $\Phi_{f_\a}(M)\vert_{X_\a}$ and as $\ba0$-module the space 
 $(\wE M)_\a$ is equivalent to $\Phi_{f_\a}(M)\vert_{X_\a}$ as
 $\pi_1(\Xa\setminus X_\a)$-module, the operators $A_{\a,\b}$ and 
 $A_{\b,\a}$, $\a\to\b$ are proportional to canonical maps $can$ and 
 $var$  between $D_{X_\a}$-modules  $\Phi_{f_\a}(M)\vert_{X_\a}$ and
 $\Psi_{f_\a}(M)\vert_{X_\a}$. Thus, the crucial calculations are the
 calculation of the vanishing cycles functor 
 $\Psi_{f_\a}(M)\vert_{X_\a}$ with the action of monodromy operator on it.
  These calculations were made in \cite{Kh},  under 
 unipotent monodromy assumption, that is, for zero weights $\l_i$.
 
 As a result of a more general calculation of the specialization functors, 
 we show in the next section that the vanishing cycles calculations 
 remain valid also under a weaker non-resonance assumption.  This 
 allows to establish an equivalence of categories for 
 non-resonance weighted $D$-modules as well.
   
%%%%%%%%%%%%%%%%%%%%%%%%%%%%%%%%%%%%%%%%%%%%%%
%%%%%%%%%%%%%%%%%%%%%%%

\bigskip

\setcounter{equation}{0}

\section{Specialization of quivers and $D$-modules}

\bigskip

{\large\bf 5.1. Specialization of quivers.}

\bigskip

Let $X_\a$ be a stratum of stratification $X_\A$. Its closure $\Xa$ is an 
affine subspace of $X=\C^N$, and the total space of normal bundle  $T_{\Xa}X$ 
 inherits the affine structure and can be naturally identified with a product
 of the affine space $\Xa$ and the linear space $X/\Xa$: $T_{\Xa}X\approx
\Xa\times X/\Xa$. 
 Both $\Xa$ and $X/\Xa$ are stratified spaces: the strata of $\Xa$ are
  $X_\b, X_\b\subset \Xa$, 
$\b\in\Ga$, and the strata of $X/\Xa$ are the
 projections $p_\a(X_\b)$, $\b\in\Ga$, $X_\b\cap\Xa\not= 0$. Here $p_\a:
 \C^N=X\to X/\Xa$ is a natural projection. Let us denote
  strata of $\Xa$ as $X_{\b'}$ and strata of $X/\Xa$ as $X_{\b''}$. 
 Note that the stratification of $X/\Xa$ is linear: the closure of any 
 stratum contains zero vector.

The stratification $T_{\Xa}\A$ of the space $T_{\Xa}X$ identified as 
 affine space with $\Xa\times X/\Xa$ is described by certain direct
 products of strata of $\Xa$ and $X/\Xa$:
 $(X_{\b'}\times X_{\b''})\in T_{\Xa}\A$ if there exists such a stratum
 $X_\b$ of stratification $\A$ for $\C^N$ that
 $$X_{\b'}=\Xa\cap X_\b , \qquad X_{\b''}= p_\a(X_\b)\ .$$
 One can see that  the stratification $T_{\Xa}\A$ is also described 
 by intersections of hyperplanes $\tilde{H}_i, i \in I: 
 \tilde{H}_i = (H_i\cap \Xa)\times p_\a(H_i)$.
The passage from stratification $\A$ to $T_{\Xa}\A$ is illustrated
 by a two-dimensional example below.

\medskip
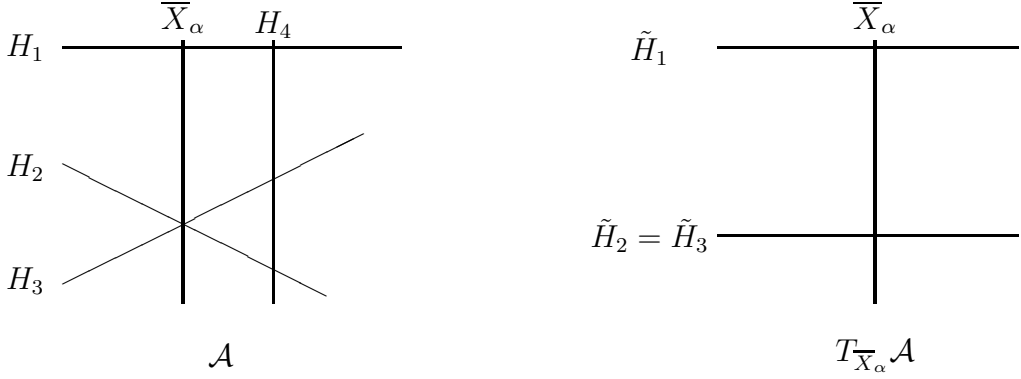
\unitlength 1mm
\linethickness{0.4pt}

\begin{picture}(120.00,45.00)
\put(07.00,39.00){\line(1,0){45}}
\put(94.00,39.00){\line(1,0){40}}
\put(94.00,14.00){\line(1,0){40}}
\put(22.90,40.00){\line(0,-1){35}}
\put(23.00,40.00){\line(0,-1){35}}
\put(23.10,40.00){\line(0,-1){35}}
\put(35.00,40.00){\line(0,-1){35}}
\put(07.00,23.50){\line(2,-1){35}}
\put(07.00,7.50){\line(2,1){40}}
\put(114.90,40.00){\line(0,-1){35}}
\put(115.00,40.00){\line(0,-1){35}}
\put(115.10,40.00){\line(0,-1){35}}
\put(23.00,43.00){\makebox(0,0)[cc]{$\overline{X}_\alpha$}}
\put(02.00,39.00){\makebox(0,0)[cc]{$H_1$}}
\put(85.00,39.00){\makebox(0,0)[cc]{$\tilde{H}_1$}}
\put(02.00,23.00){\makebox(0,0)[cc]{$H_2$}}
\put(02.00,8.00){\makebox(0,0)[cc]{$H_3$}}
\put(35.00,42.00){\makebox(0,0)[cc]{$H_4$}}
\put(85.00,14.00){\makebox(0,0)[cc]{$\tilde{H}_2=\tilde{H}_3$}}
\put(115.00,43.00){\makebox(0,0)[cc]{$\overline{X}_\alpha$}}
\put(28.00,-02.00){\makebox(0,0)[cc]{${\cal A}$}}
\put(115.00,-02.00){\makebox(0,0)[cc]{$T_{\overline{X}_\alpha}{\cal A}$}}

\end{picture}
\bigskip
%\begin{center}
%Fig.1
%\end{center}

 Correspondingly, the graph $\Gamma_{T_{\Xa}\A}$ is a subgraph of direct
 product of the graph $\Gamma'_\A$ of the stratification of the space
 $\Xa$ and of the graph $\Gamma''\A$ of the stratification of the space
 $X/\Xa$: $\Gamma_{T_{\Xa}\A} \in \Gamma'_\A\times \Gamma''_\A$. It means
 that the vertices of  $\Gamma_{T_{\Xa}\A}$ are numerated by some pairs
 $(\b',\b'')$ of the vertices of  $\Gamma'\A$ and of $\Gamma''_\A$, and
 the arrows could be of two possible types:
$$(\b',\b'')\to(\g',\b'')\ ,\qquad {\rm if}\;\; 
\b'\to\g'\in\ {\rm arr}(\Gamma'_\A)\ ,$$
 $$(\b',\b'')\to(\b',\g'')\ ,\qquad {\rm if}\;\; 
\b''\to\g''\in\ {\rm arr}(\Gamma''_\A)\ .$$
 Since the stratification $T_{\Xa}\A$ is given by intersections
 of hyperplanes $\tilde{H}_i$, we have by definition 1.1 the 
 quiver category $Qui(T_{\Xa}\A)$.

Let now $\VA =\{V_\b, A_{\b,\g}\}$ be a quiver from the category $\Qui$.
 We define a quiver $Sp_\a\VA  = 
 \{V_{(\b',\b'')}, A_{(\g',\g''),(\b',\b'')}\}$
 as follows. The space $V_{(\b',\b'')}$ is the following direct sum of 
 linear spaces $V_\b$:
\beq
V_{(\b',\b'')}=\oplus_{\b:X_\b\cap\Xa=X_{\b'}, p_\a(X_\b)=X_{\b''}}V_\b\ .
\label{4.3}
\eeq
 The operators $A_{(\b',\b''),(\g',\g'')}$ coincide with operators from
 the quiver $\VA$:
\beq
A_{(\g',\b''),(\b',\b'')}(v_\b)=
\oplus_{\g:X_\g\cap\Xa=X_{\g'}, p_\a(X_\g)=X_{\b''}}A_{\g,\b}(v_\b)\ ,
\label{4.4}
\eeq
\beq
A_{(\b',\g''),(\b',\b'')}(v_\b)=
\oplus_{\g:X_\g\cap\Xa=X_{\b'}, p_\a(X_\g)=X_{\g''}}A_{\g,\b}(v_\b)\ ,
\label{4.5}
\eeq
where the vector $v_\b$ belongs to linear space $V_\b$ such that
 $X_\b\cap\Xa=X_{\b'}, p_\a(X_\b)=X_{\b''}$.
 A map $Sp_\a$ is also defined on the weights. For a collection of weights
 $\{\l\} =\l_i =\l(H_i)$ of stratification $\A$ the collection
 $\{ Sp_\a(\l)\}$ is defined by the rule $Sp_\a(\l)_i =\sum_{j:\tilde{H}_j=
\tilde{H}_i}\l(H_i)$.
 
\begin{proposition}

{\em (i)} The relations \rf{4.3}--\rf{4.5} define an exact functor
$$Sp_\a\ : Qui(\A)\to Qui(T_{\Xa}\A)$$
 which commutes with the duality ${}^t$;

{\em (ii)} For any collection $\l$ of weights for the configuration $\A$, 
  \beq
 Sp_\a(M_\l) =M_{Sp_\a(\l)}\ .
 \label{4.6}
\eeq   
\end{proposition}

{\bf Proof.} In order to check {\it (i)} it is suficient to check that
 the operators defined by the relations \rf{4.4} and \rf{4.5} satisfy 
 the relations \rf{1.2} for the configuration $T_{\Xa}\A$. These are
 linear algebra calculations. Analogously,
 for the proof of {\it (ii)} one has to check the formulas \rf{2.2.10},
 \rf{2.2.11} or to use the functorial definition of the modules
 $M_\l$ and $M_{Sp_\a(\l)}$.

Let now $\oX_{\a_1}\subset\oX_{\a_2}\subset X=\C^N$  be a flag of 
 inclusions of closed strata to $X$. Then the spaces 
 $T_{\oX_{\a_1}}(T_{\oX_{\a_2}})$ and $T_{T_{\oX_{\a_1}}\oX_{\a_2}}
(T_{\oX_{\a_1}}X)$ are canonically isomorphic to a direct product
 $\oX_{\a_1}\times \oX_{\a_2}/\oX_{\a_1} \times X/\oX_{\a_2}$.
 Let us call it the space of flags of normal covectors and denote
 by $T_{\oX_{\a_2}\subset \oX_{\a_2}}$. It is a stratified space
 with a stratification analogous to described above. We denote it by
$T_{\oX_{\a_2}\subset \oX_{\a_2}}\A$. One can easily see that the compositions
$$Sp_{\a_1}Sp_{\a_2} : Qui(\A)\to Qui(T_{\oX_{\a_2}\subset \oX_{\a_2}}\A)$$
 and  
$$Sp_{\a_2}Sp_{\a_1} : Qui(\A)\to Qui(T_{\oX_{\a_2}\subset \oX_{\a_2}}\A)\ ,$$
    are canonically isomorphic. Here the functor $Sp_{\a_2}$ in the second
 composition is the functor of specialization to the stratum 
$\oX_{\a_1}\times \oX_{\a_2}/\oX_{\a_1}$. 

Analogously, for a three-flag 
$\oX_{\a_1}\subset\oX_{\a_2}\subset \oX_{\a_3}\subset X=\C^N$ the morphism
 of associativity is identical on the linear spaces of a quiver, so  
 for any flag
$$\oX_{\a_1}\subset\oX_{\a_2}\subset\ldots\subset \oX_{\a_n}\subset X=\C^N$$
 of closed strata we have a specialization functor 
$$Sp_{\a_1,\ldots \a_n}: Qui(\A)\to Qui(T_{\oX_{\a_1}\subset\ldots
\subset\oX_{\a_n}})$$
 which is described by itereted relations \rf{4.3}--\rf{4.5} and satisfy
the properties mentioned in Proposition 4.1.
%%%%%%%%%%%%%%%%%%%%%%%%%%%%%%%%%%%%
%%%%%%%%%%%%%%%%%%%%%%%%%%%%%%%%%%%%%
\bigskip

{\large\bf 5.2. Specialization of quiver $D$-modules.}

\bigskip

Let us first consider the case of {\it linear} stratification in linear space 
 $\C^N$, that is we suppose that all hyperplanes $H_i$ contain an origin.
 Let $I_\a$ be an ideal of $\Xa$: $f \in I_\a \Leftrightarrow f|_{\Xa}=0$
 and let $F^\bullet_\a(D_X)$ be a filtration of the ring $D_X$ corresponding
 to ideal $I_\a$ \cite{Ksh}:
$$F^k_\a(D_X)=\{ P\in D_X: P(I_\a^j\subset I_\a^{j+k} 
 {\rm for} \: {\rm any}\: j\}\ .$$
Here $I_\a^k=\C[X]$ for $k\leq 0$.
 Then $F^k_\a(D_X)/F^{k+1}_\a(D_X)$ is isomorphic to a space of differential
 operators over $T_{\Xa}X$ homogenuous of order $k$ and the ring 
 $Gr\ F_\a(D_X)$ is isomorphic to the ring $D_{T_{\Xa}X}$. 

Let now $\VA$ be a quiver from $Qui(\A)$ and $E\VA$ be corresponding 
 $D_X$-module. We define a filtration $F^\cdot_\a E\VA$ of $E\VA$ in the 
 following way:
\beq
F_\a^nE\VA=\sum_{k=0}^N\sum_{\b: codim_{X_\b} X_\b\cap \Xa =k}
F_\a^{n-k}(D_X)\cdot(V_\b\ot \O_\b)\ .
\label{4.7}
\eeq
One can check that the filtration $F^\bullet_\a E\VA$ is good, that is 
 satisfy the following conditions:
\smallskip

 (i) $F_\a^k(D_X)\cdot F_\a^jE\VA\subset F_\a^{k+j}E\VA$;
\smallskip

(ii) $F_\a^k(D_X)\cdot F_\a^jE\VA\ = F_\a^{k+j}E\VA$ for $j,k\geq 0$ and for
 $j\leq -N, k\leq 0$;
\smallskip

(iii) $\cup_j F_\a^jE\VA = E\VA$;
\smallskip

(iv) $F_\a^jE\VA$ is coherent $F^0_\a(D_X)$-module.
\smallskip

The following proposition describes $Gr\ F_\a(E\VA)$ 
as $D_{T_{\Xa}X}$- module. It follows directly from the definition
 of quiver's $D$-module $E\VA$ and from the definition of the filtration.
  
\begin{proposition}
There is a canonical isomorphism of $D_{T_{\Xa}X}$-modules:
\beq
Gr\ F_\a(E\VA)\approx ESp_\a\ \VA\ .
\label{4.8}
\eeq
\end{proposition}

 As it follows from the relations \rf{4.3}--\rf{4.5}, the isomorphism 
 \rf{4.8} means that $D_{T_{\Xa}X}$- module $Gr\ F_\a(E\VA)$ is generated
 by vectors $v_\a\ot \o_{\b'}\wedge\o_{\b''}$ where $\b\in\ver (\Ga)$, 
 $v_\b \in V_\b$, $\o_{\b'}$ is a top form for $X_\b\cap \Xa$, $\o_{\b''}$
 is a top form for $X_\b/X_\b\cap\Xa = p_\a(X_\b)$   with the following 
 relations:
\beq
L_{\b'}(v_\b\ot\o_{\b'}\wedge\o_{\b''})=
\sum_{{\g: \b\to\g \atop X_\g\cap\Xa\not=X_\b\cap\Xa}}
\frac{L_{\b'}(\ho_{\g'})\wedge\o_{\b'}}{\ho_{\g'}\wedge\o_{\g'}}
A_{\g,\b}(v_b)\ot\o_{\g'}\wedge\o_{\b''}
\label{4.9}
\eeq
where $L_{\b'}$ is linear vectorfied along stratum $X_\b\cap\Xa$,
 $\ho_{\g'}$ is a top form for $\Xa/\oX_\g\cap\Xa$, $\o_{\g'}$ is a
 top form for $\oX_\g\cap\Xa$;
\beq
L_{\b''}(v_\b\ot\o_{\b'}\wedge\o_{\b''})=
\sum_{{\g: \b\to\g \atop X_\g\cap\Xa =X_\b\cap\Xa}}
\frac{L_{\b''}(\ho_{\g''})\wedge\o_{\b''}}{\ho_{\g''}\wedge\o_{\g''}}
A_{\g,\b}(v_b)\ot\o_{\b'}\wedge\o_{\g''}
\label{4.10}
\eeq
where $L_{\b''}$ is linear vectorfied in $X/\Xa$ along  $X_\b/X_\b\cap\Xa$,
 $\ho_{\g''}$ is a top form for $(X/\Xa)/(X_\g/\oX_\g\cap\Xa)$, 
 $\o_{\g''}$ is a
 top form for $X_\g/\oX_\g\cap\Xa$;
\beq
f_{\b'}(v_\b\ot\o_{\b'}\wedge\o_{\b'')}=
\sum_{{\g: \g\to\b \atop X_\g\cap\Xa\not=X_\b\cap\Xa}}
\frac{df_{\b'}\wedge\ho_{\g'}\wedge\o_{\b'}}{\ho_{\g'}\wedge\o_{\g'}}
A_{\g,\b}(v_b)\ot\o_{\g'}\wedge\o_{\b''}
\label{4.11}
\eeq
where $f_{\b'}$ is linear function on $\Xa$ equal to zero on $X_\b\cap\Xa$,
 with the same $\ho_{\g'}$,  $\o_{\g'}$ as for \rf{4.9};
\beq
f_{\b''}(v_\b\ot\o_{\b'}\wedge\o_{\b''})=
\sum_{{\g: \g\to\b \atop X_\g\cap\Xa =X_\b\cap\Xa}}
\frac{df_{\b''}\wedge\ho_{\g''}\wedge\o_{\b''}}{\ho_{\g''}\wedge\o_{\g''}}
A_{\g,\b}(v_b)\ot\o_{\b'}\wedge\o_{\g''}
\label{4.12}
\eeq
where $f_{\b''}$ is linear function on $X/\Xa$ 
equal to zero on 
  $X_\b/X_\b\cap\Xa$ with the same 
 $\ho_{\g''}$ and  
 $\o_{\g''}$ as in \rf{4.10}.
\smallskip

  We can drop the assumption of linearity for the stratification if
 we pass first to analytic $D$-modules, and thus kill the strata which
 do not intersext with $\Xa$, and then use the well known equivalence between
 analytic and algebraic $D$-modules. Alternatively, we could first localize
 the initial $D$-module to a Zariski open set which does not contain strata
 not intersecting with $\Xa$, and then use the construction above.

In order to compare $D_{T_{\Xa}X}$- module $Gr\ F_\a(E\VA)$ with the  
 Malgrange-Kashiwara specialization of the $D_X$-module $E\VA$, we should 
 investigate the eigenvalue of the target vectorfiled $\theta^\a$ which 
 satisfy the relations $\theta^\a|_{I_\a/I_\a^2}=id$ \cite{Ksh}.
 In local coordinates $\{ \overline{z}^\a,\overline{w}^\a\} = 
 \{ z_1^\a,\ldots , z^\a_{dim\ X_\a}, w_1^\a, \ldots , w_{codim\ X_\a}^\a\}$
 satisfying \rf{2.3} it looks like $\theta^\a =\sum_j w_j^\a
\frac{\partial}{\partial w_j^\a}$. Let also $\l=\{ \l_1,\ldots , \l_n\}$
 be a collection of weights and $M_\l$ be a quiver constructed in section 3.1
 (see \rf{2.2.5}).
\begin{proposition}
\beq
\theta^\a|{Gr^kF_\a EM_\l} =\l_\a+k \ .
\label{4.13}
\eeq
\end{proposition}
The proof follows from the relations \rf{4.9}, \rf{4.10}.

Correspondingly, for $D_X$-module $E(M_{\l,\b})$ supported on the
 stratum $\oX_\b$ ($\oX_\b\cap\Xa\not= 0$), we have
$$  \theta^\a|{Gr^kF_\a EM_{\l,\a}} =\l_{\b',\b}+k \ ,
$$
where $X_{\b'}=X_\b\cap\Xa$.
\begin{corollary}

{\em (i)} Let the number $\l_\a=\sum_{i:H_i\supset X_\a} 
\not\in {\bf Z}\setminus 0$.
Then
\beq
Sp_{\Xa}E(M_\l)\approx E\  Sp_\a(M_\l)\ ;
\label{4.14}
\eeq

{\em (ii)} Let the  nonresonance condition 
 takes place: $\l_{\b,\g}-\l_{\nu,\varepsilon}\not\in{\bf Z}\setminus 0$.
Then
\beq
Sp_{\Xa}E(M)\approx E\  Sp_\a(M)
\label{4.15}
\eeq
for any $M\in\Qui$.
\end{corollary}
Here  $Sp_\Xa$ denotes the Malgrange-Kashiwara specialization functor, 
 \cite{Ksh}.

\bigskip

 {\bf Proof of Theorem 4.1}. Following the scheme of \cite{Kh}, we
 choose a collection of generic linear functions $f_\a$, $f_\a|_{\Xa}=0$
 and perform the glueing procedure starting from a local system on the
 open stratum $X_\ee$ and adding step by step strata of bigger codimension.
 In order to glue the stratum $X_\a$ we add 
 hyperplane $f_\a=0$ to the stratification with zero weight and with
 a condition that the support of $D$-module intersects with hyperplane
 $f_\a=0$ only by stratum $\Xa$. Part {\it (ii)} of the Corollary 4.1
 allows to calculate the linear algebra data    
 $$ \Psi_{f_\a}(M){\leftarrow\atop\rightarrow}\Phi_{f_\a}(M)$$
 in terms of Malgrange-Kashiwara speciaization functor assuming we
 start from quiver's $D$-modules with nonresonance condition.
  
Corollary 4.1 guarantees that they coincide with the linear algebra
 data calculated by means of specialization functors in category
 $\quila$ (and thus coincide with those calculated by Beilinson's
 technique for the unipotent case in \cite{Kh}). Then we can refer to
 \cite{Kh} since we drop unipotent restriction. Alternatively, one can get an  
 independent proof after the formulation of an analog of the  
 Beilinson-Kashiwara
 equivalence theorem for the category $\Qui$ and functors $Sp_\a$ from
 section 4.1, which is purely combinatorial.

S.K.: Institute for Theorerical and Experimental Physics, Cheremushkinskaya 25, 
117259 Moscow, Russia;\ khor@heron.itep.ru

V.S.: Max-Planck-Institut f\"ur Mathematik, Gottfried-Claren-Stra\ss e 26, 
53225 Bonn, Germany;\ vadik@mpim-bonn.mpg.de

\end{document}